\documentclass[12pt]{article}
\usepackage[letterpaper, left=1.5in, right=1in,top=1in,bottom=1in]{geometry}
\usepackage{tikz}
\usepackage{chngcntr}
\usepackage{authblk}
\counterwithin{figure}{section}
\usetikzlibrary{arrows}
\usepackage{placeins}	
\usepackage{amsfonts}
\usepackage{graphicx}
\usepackage{standalone}
\usepackage{amsmath, amsthm, amssymb}
\usepackage[capitalize, noabbrev]{cleveref}
\newtheorem{thrm}{Theorem}[section]

\newtheorem{cor}[thrm]{Corollary}

\newenvironment{pf}           {\noindent{\bf Proof:} }%
                                {\null\hfill$\Box$\par\medskip}
                           
\usepackage[nottoc]{tocbibind}

\begin{document}

\title{Orthogonal Colourings of Tensor Graphs}
\author{Kyle MacKeigan}
\affil{Department of Mathematics, Dalhousie University}
\affil{Email: \textit{Kyle.m.mackeigan@gmail.com}}
\maketitle

\begin{abstract}
In this paper, perfect k-orthogonal colourings of tensor graphs are studied. First, the problem of determining if a given graph has a perfect 2-orthogonal colouring is reformulated as a tensor subgraph problem. Then, it is shown that if two graphs have a perfect $k$-orthogonal colouring, then so does their tensor graph. This provides an upper bound on the $k$-orthogonal chromatic number for general tensor graphs. Lastly, two other conditions for a tensor graph to have a perfect $k$-orthogonal colouring are given.
\end{abstract}

\section{Introduction}

Two colourings of a graph are \textit{orthogonal} if when two elements are coloured with the same colour in one of the colourings, then those elements receive distinct colours in the other colouring. Archdeacon, Dinitz, and Harary \cite{archdeacon1985orthogonal} originally studied this type of colouring, in the context of edge colourings. Then, Caro and Yuster \cite{caro1999orthogonal} revisited this concept, this time in the context of vertex colouring. In this paper, the vertex variation is studied. 

A $k$-\textit{orthogonal colouring} of a graph $G$ is a collection of $k$ mutually orthogonal vertex colourings. For simplicity, a $2$-orthogonal colouring is called an \textit{orthogonal colouring}. The \textit{$k$-orthogonal chromatic number} of a graph $G$, denoted by $O\chi_k(G)$, is the minimum number of colours required for a proper $k$-orthogonal colouring. Again for simplicity, the $2$-orthogonal chromatic number is simply denoted by $O\chi(G)$ and simply called the \textit{orthogonal chromatic number}.

For a graph $G$ with $n$ vertices, $O\chi_k(G)\geq \lceil\sqrt{n}\,\rceil$. If $G$ has $n^2$ vertices and $O\chi_k(G)=\lceil\sqrt{n^2}\,\rceil=n$, then $G$ is said to have a \textit{perfect $k$-orthogonal colouring}. A perfect $2$-orthogonal colouring is simply called a \textit{perfect orthogonal colouring}. Perfect orthogonal colourings are of particular importance because they have applications to independent coverings \cite{mackeigan2021independent} and scoring games \cite{andres2019orthogonal}. 

Due to these applications, research is focused on determining which graphs have perfect $k$-orthogonal colourings. For instance, Caro and Yuster \cite{caro1999orthogonal} constructed graphs having perfect $k$-orthogonal colourings by using orthogonal Latin squares. Ballif \cite{ballif2013upper} studied upper bounds on sets of orthogonal colourings. Whereas Janssen and the author \cite{janssen2020orthogonal} studied perfect $k$-orthogonal colourings of circulant graphs. 

In this paper, the tensor graph product is used to construct graphs having perfect $k$-orthogonal colourings. The \textit{tensor product} of two graphs $G$ and $H$, denoted by $G\times H$, has vertex set $V(G)\times V(H)$, and two vertices $(u_1,v_1)$ and $(u_2,v_2)$ in $G\times H$ are adjacent if and only if $u_1u_2\in E(G)$ and $v_1v_2\in E(H)$. If $G$ is a graph that is created by the tensor product of two graphs, then $G$ is called a \textit{tensor graph}.

For graph products, there can be a clear relationship between colourings  of the component graphs and colourings of the product graph. For example, it is well known that the chromatic number of a Cartesian graph (see Section 3 for a formal definition) is the minimum of the chromatic numbers of the components. For tensor graphs, it was conjectured by Hedetniemi \cite{hedetniemi1966homomorphisms} that the same result would be true. However, Hedetniemi's conjecture was recently disproved by Shitov \cite{shitov2019counterexamples}.

In this paper, it is shown that for a graph $G$, $O\chi(G)\leq n$ if and only if $G$ is a subgraph of the tensor graph $K_n\times K_n$. Then, it is shown that if two graphs have a perfect $k$-orthogonal colouring, then so does their tensor graph. For $k=2$, only one component is required to have a perfect orthogonal colouring. Similarly, it is shown that this condition can be relaxed for perfect $k$-orthogonal colourings. These results were first obtained in the authors doctoral dissertation \cite{mackeigan2021exploration}.

\section*{Acknowledgements}

The author would like to thank his supervisor Jeannette Janssen for her helpful comments and suggestions on this paper.

\section{Perfect Orthogonal Colourings}

In this section, perfect orthogonal colourings of tensor graphs are studied. To start, it is shown that a graph $G$ has $O\chi(G)\leq n$ if and only if is a subgraph of $K_n\times K_n$, which in this paper, is denoted by $G\subseteq K_n\times K_n$. Therefore, a graph $G$ with $m^2$ vertices has a perfect orthogonal colouring if and only if it is a subgraph of $K_m\times K_m$.

\begin{thrm}\label{Theorem: Categorization}
For a graph $G$, $O\chi(G)\leq n$ if and only if $G\subseteq K_n\times K_n$.
\end{thrm}
\begin{pf}
For $1\leq i,j\leq n$, let $(i,j)$ denote the vertices of the graph $K_n\times K_n$. First, suppose that $G\subseteq K_n\times K_n$. It will be shown that $K_n\times K_n$ has a perfect orthogonal colouring. If this is the case, then the orthogonal colouring of $K_n\times K_n$ restricted to $G$ is an orthogonal colouring of $G$ using $n$ colours, giving $O\chi(G)\leq n$.

Assign the vertex $(i,j)$ the colour $i$ in the first colouring and the colour $j$ in the second colouring. For example, this orthogonal colouring is applied to $K_3\times K_3$ in Figure \ref{Figure: Orthogonal Colouring of K3timesK3}. Displayed next to each vertex are the colours assigned in the first and second colouring. 

Note that this assignment of colours has no orthogonal conflicts. It remains to check that there are no colour conflicts. Now, by the definition of the tensor product, for $1\leq i_1,i_2,j_1,j_2\leq n$, two vertices $(i_1,j_1)$ and $(i_2,j_2)$ in $K_n\times K_n$ are adjacent if and only if $i_1\neq i_2$ and $j_1\neq j_2$. Therefore, there are also no colour conflicts.

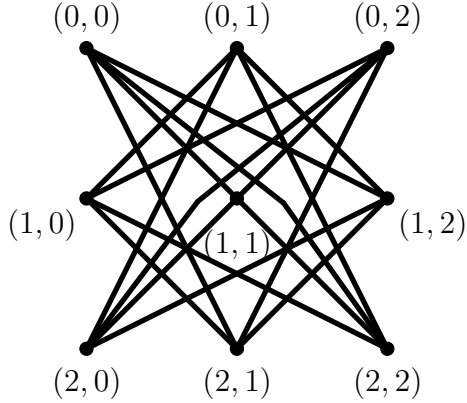
\begin{figure}[h!]
\centering
\begin{tikzpicture}[line cap=round,line join=round,>=triangle 45,x=1.0cm,y=1.0cm]
\draw [line width=2.pt] (0.,4.)-- (2.,0.);
\draw [line width=2.pt] (0.,2.)-- (2.,4.);
\draw [line width=2.pt] (0.,2.)-- (2.,0.);
\draw [line width=2.pt] (2.,4.)-- (0.,0.);
\draw [line width=2.pt] (4.,4.)-- (2.,0.);
\draw [line width=2.pt] (4.,2.)-- (2.,4.);
\draw [line width=2.pt] (4.,2.)-- (2.,0.);
\draw [line width=2.pt] (4.,0.)-- (2.,4.);
\draw [line width=2.pt] (0.,4.)-- (4.,2.);
\draw [line width=2.pt] (4.,4.)-- (0.,2.);
\draw [line width=2.pt] (4.,0.)-- (0.,2.);
\draw [line width=2.pt] (0.,0.)-- (4.,2.);
\draw [line width=2.pt] (0.,4.)-- (4.,0.);
\draw [line width=2.pt] (4.,4.)-- (0.,0.);
\draw [line width=2.pt] (0.,4.)-- (2.62,1.94);
\draw [line width=2.pt] (4.,0.)-- (2.62,1.94);
\draw [line width=2.pt] (4.,4.)-- (1.48,1.96);
\draw [line width=2.pt] (1.48,1.96)-- (0.,0.);
\draw (-0.6,4.75) node[anchor=north west] {$(0,0)$};
\draw (1.4,4.75) node[anchor=north west] {$(0,1)$};
\draw (3.4,4.75) node[anchor=north west] {$(0,2)$};
\draw (-1.2,2) node[anchor=north west] {$(1,0)$};
\draw (1.4,1.74) node[anchor=north west] {$(1,1)$};
\draw (4,2) node[anchor=north west] {$(1,2)$};
\draw (-0.6,-0.1) node[anchor=north west] {$(2,0)$};
\draw (1.4,-0.1) node[anchor=north west] {$(2,1)$};
\draw (3.4,-0.1) node[anchor=north west] {$(2,2)$};
\begin{scriptsize}
\draw [fill=black] (0.,0.) circle (2.5pt);
\draw [fill=black] (0.,2.) circle (2.5pt);
\draw [fill=black] (0.,4.) circle (2.5pt);
\draw [fill=black] (2.,0.) circle (2.5pt);
\draw [fill=black] (2.,2.) circle (2.5pt);
\draw [fill=black] (2.,4.) circle (2.5pt);
\draw [fill=black] (4.,4.) circle (2.5pt);
\draw [fill=black] (4.,2.) circle (2.5pt);
\draw [fill=black] (4.,0.) circle (2.5pt);
\end{scriptsize}
\end{tikzpicture}
\caption{Orthogonal Colouring of $K_3\times K_3$}
\label{Figure: Orthogonal Colouring of K3timesK3}
\end{figure}

Now, suppose that $O\chi(G)\leq n$ and that $(g_1,g_2)$ is an orthogonal colouring of $G$ using the colours $\{1,2,\dots,n\}$. To show that $G\subseteq K_n\times K_n$, an injective map that preserves edges is required. Let $F:G\to K_n\times K_n$ by $F(v)=(g_1(v),g_2(v))$. It is now shown that $F$ is injective and preserves edges.

 Since $(g_1,g_2)$ is an orthogonal colouring of $G$, each colour pair is only assigned once. Thus, $F$ is injective. Now, if $v_1v_2\in E(G)$, then $g_1(v_1)\neq g_1(v_2)$ and $g_2(v_1)\neq g_2(v_2)$ because $g_1$ and $g_2$ are proper. Therefore, $(g_1(v_1),g_2(v_1))(g_1(v_2),g_2(v_2))\in E(K_n\times K_n)$ by the definition of the edges in $K_n\times K_n$. Thus, $F$ preserves edges. Since $F$ is injective and preserves edges, $G\subseteq K_n\times K_n$.
\end{pf}

Theorem \ref{Theorem: Categorization} gives a way to reformulate the problem of determining if a graph has a perfect orthogonal colouring. This will be used later with the following theorem to obtain an upper bound on the orthogonal chromatic number of general tensor graphs. The following theorem shows that if one component has a perfect orthogonal colouring and the other has a square number of vertices, then their tensor graph has a perfect orthogonal colouring.

\begin{thrm}\label{Theorem: Tensor Product}
If $G$ has $n^2$ vertices, $H$ has $m^2$ vertices, and $O\chi(G)=n$, then $O\chi(G\times H)=nm$.
\end{thrm}
\begin{pf}
Label $V(G)=\{v_k:0\leq k <n^2\}$ and $V(H)=\{(u_i,u_j):0\leq i,j<m \}$. Let $f=(f_1,f_2)$ be a proper orthogonal colouring of $G$ where $f_1$ and $f_2$ use the colours $\{0,1,\dots, n-1\}$. It is shown that $g=(g_1,g_2)$ is an orthogonal colouring of $G\times H$ using $nm$ colours, where:
\begin{align*}
g_1((v_k,(u_i,u_j)))&=f_1(v_k)+in\\
\text{and}\\
g_2((v_k,(u_i,u_j)))&=f_2(v_k)+jn.
\end{align*}
First, it is shown that $g$ has no orthogonal conflicts. Let $v_{k_1},v_{k_2}\in V(G)$ and let $(u_{i_1},u_{j_1}),(u_{i_2},u_{j_2})\in V(H)$. If $g((v_{k_1},(u_{i_1},u_{j_1})))=g((v_{k_2},(u_{i_2},u_{j_2})))$, then:
\begin{align}
f_1(v_{k_1})+i_1n&=f_1(v_{k_2})+i_2n \label{Equation: ONE}\\
\nonumber\text{and}\\
f_2(v_{k_1})+j_1n&=f_2(v_{k_2})+j_2n. \label{Equation: TWO}
\end{align}
Without loss of generality, suppose that $i_1<i_2$. Then it follows that:
\begin{align*}
f_1(v_{k_1})+i_1n &<n+i_1n\\
&\leq i_2 n\\
&\leq f_1(v_{k_2})+i_2n.
\end{align*}
Therefore, $f_1(v_{k_1})+i_1n<f_1(v_{k_2})+i_2n$, which contradicts Equation \eqref{Equation: ONE}, thus $i_1=i_2$. A similar argument shows that $j_1=j_2$. Substituting $i_1=i_2$ and $j_1=j_2$ into Equations \eqref{Equation: ONE} and \eqref{Equation: TWO}, gives $f_1(v_{k_1})=f_1(v_{k_2})$ and $f_2(v_{k_1})=f_2(v_{k_2})$. Hence, $v_{k_1}=v_{k_2}$ because $f$ is an orthogonal colouring of $G$. Thus, $(v_{k_1},(u_{i_1},u_{j_1}))=(v_{k_2},(u_{i_2},u_{j_2}))$.

It remains to show that $g_1$ and $g_2$ are proper colourings of $G\times H$. Suppose that $v_{k_1}v_{k_2}\in E(G)$ and $(u_{i_1},u_{j_1})(u_{i_2},u_{j_2})\in E(H)$. If $i_1=i_2=i$, then since $f_1$ is a proper colouring of $G$, $g_1((v_{k_1},(u_{i_1},u_{j_1})))= f_1(v_{k_1})+in\neq f_1(v_{k_2})+in=g_1((v_{k_2},(u_{i_2},u_{j_2})))$. Thus, there are no colour conflicts between these verices.

Now, without loss of generality, suppose that $i_1<i_2$. Then it follows that $g_1((v_{k_1},(u_{i_1},u_{j_1})))=f_1(v_{k_1})+i_1n<n+i_1n\leq i_2n\leq f_1(v_{k_2})+i_2n = g_1((v_{k_1},(u_{i_2},u_{j_2})))$. Hence, $g_1((v_{k_1},(u_{i_1},u_{j_1})))<g_1((v_{k_1},(u_{i_2},u_{j_2})))$. Thus, there are no colour conflicts between these vertices. Therefore, $g_1$ is a proper colouring. A similar argument shows that $g_2$ is proper. Thus, $g$ is an orthogonal colouring of $G\times H$. Since $G\times H$ has $n^2m^2$ vertices and $g$ uses $nm$ colours, $O\chi(G\times H)=nm$.
\end{pf}

Theorem \ref{Theorem: Tensor Product} provides a method for constructing perfect orthogonal colourings out of graphs that have perfect orthogonal colourings. On the other hand, Theorem \ref{Theorem: Categorization} gives that $K_n\times K_n$ is the maximum graph with $n$ as its orthogonal chromatic number. Combining these two results provides an upper bound on the orthogonal chromatic number of general tensor graphs.

\begin{cor}\label{Corollary: Upper Bound}
If $O\chi(G)=n$ and $O\chi(H)=m$, then $O\chi(G\times H)\leq nm$.
\end{cor}
\begin{pf}
Since $O\chi(G)=n$ and $O\chi(H)=m$, $G\subseteq K_n\times K_n$ and $H\subseteq K_m\times K_m$ by Theorem \ref{Theorem: Categorization}. Therefore, $G\times H\subseteq (K_n\times K_n)\times (K_m\times K_m)$. Since $|V(K_n\times K_n)|=n^2$, $|V(K_m\times K_m)|=m^2$, and $O\chi(K_n\times K_n)=n$, $O\chi((K_n\times K_n)\times (K_m\times K_m))=nm$ by Theorem \ref{Theorem: Tensor Product}. Therefore, $(K_n\times K_n)\times (K_m\times K_m)\subseteq K_{nm}\times K_{nm}$ by Theorem \ref{Theorem: Categorization}. Thus, $G\times H\subseteq K_{nm}\times K_{nm}$, and Theorem \ref{Theorem: Categorization} gives that $O\chi(G\times H)\leq nm$.
\end{pf}

Corollary \ref{Corollary: Upper Bound} gives an upper bound on the orthogonal chromatic number of tensor graphs in the case where the orthogonal chromatic numbers of the components are known. However, this upper bound can be far from the exact orthogonal chromatic number. For instance, $O\chi(K_n)=n$, so Corollary \ref{Corollary: Upper Bound} gives $O\chi(K_n\times K_n)\leq n^2$. However, by Theorem \ref{Theorem: Categorization}, $O\chi(K_n\times K_n)=n$.

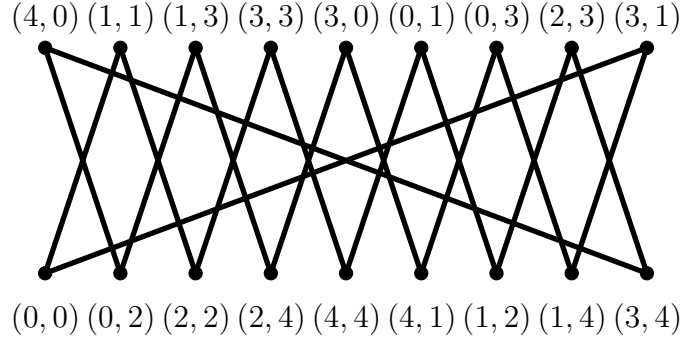
\begin{figure}[h!]
\centering
\begin{tikzpicture}[line cap=round,line join=round,>=triangle 45,x=1.0cm,y=1.0cm]
\draw [line width=2.pt] (0.,3.)-- (1.,0.);
\draw [line width=2.pt] (2.,3.)-- (1.,0.);
\draw [line width=2.pt] (2.,3.)-- (3.,0.);
\draw [line width=2.pt] (3.,0.)-- (4.,3.);
\draw [line width=2.pt] (4.,3.)-- (5.,0.);
\draw [line width=2.pt] (5.,0.)-- (6.,3.);
\draw [line width=2.pt] (6.,3.)-- (7.,0.);
\draw [line width=2.pt] (7.,0.)-- (8.,3.);
\draw [line width=2.pt] (8.,3.)-- (0.,0.);
\draw [line width=2.pt] (0.,3.)-- (8.,0.);
\draw [line width=2.pt] (0.,0.)-- (1.,3.);
\draw [line width=2.pt] (1.,3.)-- (2.,0.);
\draw [line width=2.pt] (2.,0.)-- (3.,3.);
\draw [line width=2.pt] (3.,3.)-- (4.,0.);
\draw [line width=2.pt] (4.,0.)-- (5.,3.);
\draw [line width=2.pt] (5.,3.)-- (6.,0.);
\draw [line width=2.pt] (6.,0.)-- (7.,3.);
\draw [line width=2.pt] (7.,3.)-- (8.,0.);
\draw (-.6,-.25) node[anchor=north west] {$(0,0)$};
\draw (.4,3.75) node[anchor=north west] {$(1,1)$};
\draw (1.4,-0.25) node[anchor=north west] {$(2,2)$};
\draw (2.4,3.75) node[anchor=north west] {$(3,3)$};
\draw (3.4,-0.25) node[anchor=north west] {$(4,4)$};
\draw (4.4,3.75) node[anchor=north west] {$(0,1)$};
\draw (5.4,-0.25) node[anchor=north west] {$(1,2)$};
\draw (6.4,3.75) node[anchor=north west] {$(2,3)$};
\draw (7.4,-0.25) node[anchor=north west] {$(3,4)$};
\draw (-0.6,3.75) node[anchor=north west] {$(4,0)$};
\draw (0.4,-0.25) node[anchor=north west] {$(0,2)$};
\draw (1.4,3.75) node[anchor=north west] {$(1,3)$};
\draw (2.4,-0.25) node[anchor=north west] {$(2,4)$};
\draw (3.4,3.75) node[anchor=north west] {$(3,0)$};
\draw (4.4,-0.25) node[anchor=north west] {$(4,1)$};
\draw (5.4,3.75) node[anchor=north west] {$(0,3)$};
\draw (6.4,-0.25) node[anchor=north west] {$(1,4)$};
\draw (7.4,3.75) node[anchor=north west] {$(3,1)$};
\begin{scriptsize}
\draw [fill=black] (0.,0.) circle (2.5pt);
\draw [fill=black] (1.,0.) circle (2.5pt);
\draw [fill=black] (2.,0.) circle (2.5pt);
\draw [fill=black] (3.,0.) circle (2.5pt);
\draw [fill=black] (4.,0.) circle (2.5pt);
\draw [fill=black] (5.,0.) circle (2.5pt);
\draw [fill=black] (6.,0.) circle (2.5pt);
\draw [fill=black] (7.,0.) circle (2.5pt);
\draw [fill=black] (8.,0.) circle (2.5pt);
\draw [fill=black] (0.,3.) circle (2.5pt);
\draw [fill=black] (1.,3.) circle (2.5pt);
\draw [fill=black] (2.,3.) circle (2.5pt);
\draw [fill=black] (3.,3.) circle (2.5pt);
\draw [fill=black] (4.,3.) circle (2.5pt);
\draw [fill=black] (5.,3.) circle (2.5pt);
\draw [fill=black] (6.,3.) circle (2.5pt);
\draw [fill=black] (7.,3.) circle (2.5pt);
\draw [fill=black] (8.,3.) circle (2.5pt);
\end{scriptsize}
\end{tikzpicture}
\caption{Orthogonal Colouring of $C_9\times K_2$}
\label{Figure: Orthogonal C9timesK2}
\end{figure}

On the other hand, Corollary \ref{Corollary: Upper Bound} gives a good upper bound for the bipartite double cover graphs, $G\times K_2$, where $G$ has an optimal orthogonal colouring. These graphs are of interest for other types of colourings \cite{mackeigan2020total}. For example, consider the cycle graph $C_9$ which by \cite{janssen2020orthogonal} has $O\chi(C_9)=3$. Then by Corollary \ref{Corollary: Upper Bound}, $O\chi(C_9\times K_2)\leq 6$. Which is only one off of the correct orthogonal chromatic number, illustrated in Figure \ref{Figure: Orthogonal C9timesK2}.

This concludes this sections study of perfect orthogonal colourings of tensor graphs. It remains an open problem to determine the correct orthogonal chromatic number of bipartite double cover graphs. Studying the maximal case, $K_n\times K_2$, may yield an improved upper bound.

\section{Perfect $k$-Orthogonal Colourings}

To start this section, it is shown that if two graphs have a perfect $k$-orthogonal colouring, then so does their tensor graph. The main idea behind the proof of this result is to create colour classes for the tensor graph out of the colour classes of the components. Then, the goal is to show that these new colour classes only share at most one element, and thus give an orthogonal colouring.

\begin{thrm}\label{Theorem: Both component k perfect}
If $G$ has $n^2$ vertices with $O\chi_k(G)=n$ and $H$ has $m^2$ vertices with $O\chi_k(H)=m$, then $O\chi_k(G\times H)=nm$.
\end{thrm}
\begin{pf}
For $0\leq r< k$ and $0\leq i<n$, let $G_{r,i}$ be the $i$-th colour class in the $r$-th colouring of $G$. Then, for $0\leq r<k$ and $0\leq j<m$, let $H_{r,j}$ be the $j$-th colour class in the $r$-th colouring of $H$. Next, let $I_{r,i,j}=\{(u,v) \text{ }|\text{ } u\in G_{r,i}, v\in H_{r,j}\}$. It will be shown that $C_r=\{I_{r,i,j}\text{ }|\text{ }0\leq i<n, 0\leq j<m\}$ is a collection of disjoint spanning independent sets. That is, $C_r$ is a proper colouring of $G\times H$ using $nm$ colours.

First, it is shown that each $I_{r,i,j}$ is an independent set. Let $(u_1,v_1),(u_2,v_2)\in I_{r,i,j}$. Then $u_1,u_2\in G_{r,i}$ and $v_1,v_2\in H_{r,j}$. However, $G_{r,i}$ and $H_{r,j}$ are independent sets, thus $u_1u_2\not\in E(G)$ and $v_1v_2\not\in E(H)$. Therefore, $(u_1,v_1)(u_2,v_2)\not\in E(G\times H)$. Now, let $(u,v)\in V(G\times H)$. Since $\{G_{r,i}\text{ }|\text{ }0\leq i<n\}$ is a spanning set of $G$, $u\in G_{r,i}$ for some $i$. Similarly, $v\in H_{r,j}$ for some $j$. Therefore, $(u,v)\in I_{r,i,j}$.

Now, suppose that $(u,v)\in I_{r,i_1,j_1}$ and $(u,v)\in I_{r,i_1,j_2}$. If $i_1\neq i_2$ then $u\in G_{r,i_1}$ and $u\in G_{r,i_2}$. However, this contradicts that $\{G_{r,i}\text{ }|\text{ }0\leq i<n\}$ is a colouring of $G$. Similarly, if $j_1\neq j_2$, then $v\in H_{r,j_1}$ and $v\in H_{r,j_2}$. However, this contradicts that $\{H_{r,j_1}\text{ }|\text{ }|0\leq j<m\}$ is a colouring of $H$. Therefore, there is a unique set $I_{r,i,j}$ that contains $(u,v)$. Thus, $C_r$ is a proper colouring of $G\times H$ using $nm$ colours.

It remains to show that each of the colourings are mutually orthogonal. Consider $I_{r_1,i_1,j_1}$ and $I_{r_2,i_2,j_2}$ where $r_1\neq r_2$. If $(u,v)\in I_{r_1,i_1,j_1}$ and $(u,v)\in I_{r_2,i_2,j_2}$, then $u\in G_{r_1,i_1}$ and $u\in G_{r_2,i_2}$. However, $|G_{r_1,i_1}\cap G_{r_2,i_2}|=1$, so let $u$ be this unique vertex. Similarly, $v\in H_{r_1,j_1}$ and $v\in H_{r_2,j_2}$. However, $|H_{r_1,j_1}\cap H_{r_2,j_2}|=1$, so let $v$ be this unique vertex. Therefore, there is a unique vertex $(u,v)$ in both $I_{r_1,i_1,j_1}$ and $I_{r_2,i_2,j_2}$. Hence, each of the $C_r$ are mutually orthogonal.
\end{pf}

Interestingly, the orthogonal colouring created in Theorem \ref{Theorem: Both component k perfect} works for Cartesian graphs as well. The Cartesian graph product of two graphs $G$ and $H$, denoted by $G\square H$, has vertex set $V(G)\times V(H)$, and two vertices $(u_1,v_1)$ and $(u_1,v_2)$ in $G\square H$ are adjacent if and only if $u_1u_2\in E(G)$ and $v_1v_2\in E(H)$. It is also a perfect $k$-orthogonal colouring for the strong product graph. The strong product of two graphs $G$ and $H$, denoted $G\boxtimes H$, vertex set $V(G)\times V(H)$ and edge set $E(G\boxtimes H)=E(G\square H)\cup E(G\times H)$.

\begin{cor}
If $G$ has $n^2$ vertices with $O\chi_k(G)=n$ and $H$ has $m^2$ vertices with $O\chi_k(H)=m$, then $O\chi_k(G\square H)=nm$ and $O\chi_k(G\boxtimes H)=nm$.
\end{cor}
\begin{pf}
Let $I_{r,i,j}$ be the same set as in Theorem \ref{Theorem: Both component k perfect}. Then, note that $I_{r,i,j}$ is an independent set in $G\square H$ and $G\boxtimes H$. Therefore, this results follows by applying the proof of Theorem \ref{Theorem: Both component k perfect}.
\end{pf}

Theorem \ref{Theorem: Both component k perfect} gives a method to construct perfect $k$-orthogonal colourings when both components have  a perfect $k$-orthogonal colouring. This is now used to find an upper bound on the $k$-orthogonal chromatic number of tensor product graphs. Recall Theorem \ref{Theorem: Categorization}, which gives a way to reformulate the problem as a subgraph question. Unlike perfect orthogonal colourings, for perfect $k$-orthogonal colourings, there are multiple graphs required to reformulate the problem. 

Caro and Yuster \cite{caro1999orthogonal} showed that a graph has a perfect $k$-orthogonal colouring if and only if it is a subgraph of a graph obtained by removing $k$ edge disjoint $K_n$-covers from $K_{n^2}$. Let $K_{n^2}[k]$ denote this family of graphs. Thus, for $k=2$, Theorem \ref{Theorem: Categorization} gives that $K_{n^2}[2]=K_n\times K_n$. Therefore, using the same argumentation as Corollary \ref{Corollary: Upper Bound}, but using this family of graphs, the following upper bound is obtained.

\begin{cor}\label{Corollary: k tensor}
If $O\chi_k(G)=n$ and $O\chi_k(H)=m$. then $O\chi_k(G\times H)\leq nm$.
\end{cor}
\begin{pf}
Suppose that $O\chi_k(G)=n$. Then $G\subseteq \bar G$ and $H\subseteq \bar H$ for some $\bar G \in K_{n^2}[k]$ and $\bar H\in K_{m^2}[k]$. Then, since $\bar G$ has $n^2$ vertices with $O\chi_k(G)=n$ and $\bar H$ has $m^2$ vertices with $O\chi_k(\bar H)=m$, $O\chi_k(\bar G\times \bar H)=nm$ by Theorem \ref{Theorem: Both component k perfect}. Therefore, since $G\times H\subseteq \bar G\times \bar H$, $O\chi_k(G\times H)\leq nm$ by restricting the $k$-orthogonal colouring.
\end{pf}

Corollary \ref{Corollary: k tensor} gives an upper bound on the $k$-orthogonal chromatic number of tensor graphs in the case where the orthogonal chromatic number of the components are known. Similar to Corollary \ref{Corollary: Upper Bound}, this gives better upper bounds to closer the $k$-orthogonal chromatic number of the components are to being perfect $k$-orthogonal chromatic numbers. To conclude this paper, the following theorem gives one more method to construct perfect $k$-orthogonal colourings of tensor graphs.

\begin{thrm}\label{Theorem: Perfect k tensor}
If $G$ has $n^2$ vertices, $H$ has $p^2$ vertices where $p$ is a prime, and $O\chi_k(G)=n$ with $k\leq p$, then $O\chi_k(G\times H)=np$.
\end{thrm}
\begin{pf}
Label $V(H)=\{(u_i,u_j):0\leq i,j<p \}$. For $0\leq r< k$ and $0\leq s<n$, let $I_{r,s}$ be the $s$-th colour class in the $r$-th colouring of $G$. Then, for $0\leq j<p$, let $\bar I_{r,s,j}=\{(v,(u_{i},u_{(ir+j)(\textrm{mod}~p)}))\text{ }|\text{ }v\in I_{r,s},0\leq i<p\}$. The goal is to show that $C_r=\{\bar I_{r,s,j}\text{ }|\text{ }0\leq s<n, 0\leq j<p\}$ is a partition of $G\times H$ into $np$ independent sets. That is, $C_r$ is a proper colouring of $G\times H$ using $np$ colours.

First, it is shown that each $\bar I_{r,s,j}$ is an independent set. Since each $I_{r,s}$ is an independent set in $G$, for each $v_1,v_2\in I_{r,s}$, $v_1v_2\not\in E(G)$. Thus by the definition of the tensor product graph, $(v_1,(u_{i},u_{(ir+j)(\textrm{mod}~p)}))(v_2,(u_{i},u_{(ir+j)(\textrm{mod}~p)}))\not\in E(G\times H)$. Therefore, each $\bar I_{r,s,j}$ is an independent set. Next, it is shown that $C_r$ is a partition of $G\times H$. 

Consider a vertex $(v,(u_x,u_y))$ in $G\times H$. Since $\{I_{r,s}\text{ }|\text{ } 0\leq s<n\}$ is a partition of $G$, $v\in I_{r,s}$ for some $s$. Now, notice that for $0\leq j<p$, $\{(u_{i},u_{(ir+j)(\textrm{mod}~p)})\text{ }|\text{ }0\leq i<p\}$ is a partition of $H$. Therefore, $(u_x,u_y)$ is in one of these sets. In particular, this occurs for $i=x$ and $j=y-r$. Therefore, there is a unique set $\bar I_{r,s,j}$ that contains $(v,(u_x,u_y))$. Thus, $C_r$ is a partition of $G\times H$ into independent sets.

Now it remains to show that each of the colourings are mutually orthogonal. That is, it remains to show for $r_1\neq r_2$, $s_1,s_2$ and $j_1,j_2$ fixed, that $|\bar I_{r_1,s_1,j_1}\cap \bar I_{r_2,s_2,j_2}|=1$. Since $|I_{r_1,s_1}\cap I_{r_2,s_2}|=1$, let $v$ be this vertex. Therefore, if it can be shown that $|\{(u_{i_1},u_{i_1r_1+j_1})|0\leq i_1<p\}\cap \{(u_{i_2},u_{i_2r_2+j_2})|0\leq i_2<p\}|=1$, then we are done. Since $r_1\neq r_2$, the only way $(u_{i_1},u_{i_1r_1+j_1})=(u_{i_2},u_{i_2r_2+j_2})$ is if $i_1=i_2$. Thus, call this $i$. Now, consider the following equation.
\begin{align*}
(ir_1+j_1)(\textrm{mod}~p) &= (ir_2+j_2)(\textrm{mod}~p)
\end{align*}

This equation simplifies to the following equation.
\begin{align*}
i(r_1-r_2)(\textrm{mod}~p)&=(j_2-j_1)(\textrm{mod}~p)
\end{align*}

Since $k\leq p$, and $r_1\neq r_2$, $r_1-r_2\neq 0(\textrm{mod}~p)$. Thus, $r_1-r_2=r$ and $j_2-j_1=j$. Since $p$ is a prime, $\mathbb{Z}_p$ has no zero divisors. Therefore, $ir(\textrm{mod}~p) =j(\textrm{mod}~p)$ has a unique solution, call this unique solution $(i,j)$. Thus, $(v,(i,j))$ is the unique element in the  $\bar I_{r_1,s_1,j_1}\cap \bar I_{r_2,s_2,j_2}$. Hence, the colourings are all mutually orthogonal. Since each of these colourings using $np$ colours, and $G\times H$ has $n^2p^2$ vertices, $O\chi(G\times H)=np$ as desired.
\end{pf}

\bibliographystyle{amsplain}
\bibliography{Orthogonal_Colourings_of_Tensor_Graphs_References}
\end{document}